\documentclass[11pt]{article}
\usepackage{times,epsfig,amssymb,latexsym,lscape}

{}
{}
{}

\title{Functions to Support Input and Output of Intervals}

\author{M.H. van Emden
        \thanks {Department of Computer Science, University of
                 Victoria, Canada}
        \and
        B. Moa
        \thanks {Department of Computer Science, University of
                 Victoria, Canada}
        \and
        S.C. Somosan
        \thanks{NewHeights Software Corporation, Victoria, Canada}
        }

\date{Research Report DCS-311-IR\\
Department of Computer Science\\
University of Victoria Canada
        }

\begin{document}
\maketitle
\pagestyle{empty}

\begin{abstract}
Interval arithmetic is hardly feasible
without directed rounding as provided,
for example, by the IEEE floating-point standard.
Equally essential for interval methods
is directed rounding for conversion
between the external decimal and internal binary numerals.
This is not provided by the standard I/O libraries.
Conversion algorithms exist that guarantee identity upon conversion
followed by its inverse.
Although it may be possible to adapt these algorithms 
for use in decimal interval I/O,
we argue that outward rounding in radix conversion
is computationally a simpler problem than guaranteeing identity.
Hence it is preferable to develop decimal interval I/O \emph{ab initio},
which is what we do in this paper.
\end{abstract}

\section{Introduction}

Interval arithmetic endows every computation with the authority of proof ---
the theorem being that the real-valued solution $x$
belongs to a set of reals $[a,b]$,
where $a$ and $b$ are IEEE-standard floating-point numbers.
Yet it can easily happen that we get as output
something obviously wrong such as
\verb+[0.33333,0.33333]+ when $x = 1/3$.
It may well be that a correct interval has been computed.
For example,
if $[a,b]$ is the narrowest
single-length IEEE-standard floating-point interval containing $1/3$,
and if we do not allow the output to be truncated,
we get $[a,b]$ as
\begin{center}
\begin{verbatim}
[0.333333313465118408203125,
 0.3333333432674407958984375],
\end{verbatim}
\end{center}
which is best written as
\begin{center}
\begin{verbatim}
0.3333333[13465118408203125,432674407958984375],
\end{verbatim}
\end{center}
a notation proposed and analyzed in \cite{vnmdn04}.
This is the unabridged version of
\begin{center}
\begin{verbatim}
[0.33333,0.33333],
\end{verbatim}
\end{center}
which is what we get with a typical default precision.
The untruncated decimal numerals
are exact representations of $[a,b]$
because every binary floating-point number
can be represented as a decimal numeral.
However, if one wants to ensure exact representation, then
we get a decimal for every bit;
not an economical representation.
Then we might as well write the bits themselves,
which gives
\begin{center}
\begin{verbatim}
0.010101010101010101010101[0,1].
\end{verbatim}
\end{center}

On output the only problem is that any reasonable choice of display precision
causes the standard output routines to shorten the numerals for the bounds to
the same result.
The improvement needed here is output that is aware of whether an upper or a
lower bound is to be displayed.

On input, however, we have a more serious problem.
Suppose we want to initialize an interval variable at $0.1$.
Initializing the lower and upper bounds with
\begin{verbatim}
          lb = ub = 0.1
\end{verbatim}
is guaranteed to generate an interval
that does \emph{not} contain $0.1$,
for the simple reason
that there is no floating-point number equal to $0.1$.
Thus the best the compiler can do is to produce
one of the floating-point numbers closest to $0.1$.
We don't know which one it is.
The table in Figure~\ref{fig:sioresults}
shows for each of $1/2, 1/3, \ldots, 1/10$
that the compiler picks
the upper bound of the narrowest interval
containing the fraction concerned.
To remind us that we cannot count on this to happen,
the other choice is made for $1/11$.
Thus, for input we need an algorithm that produces,
for every fraction or numeral as input,
a narrow interval containing it;
ideally the narrowest.

Floating-point numbers are normalized as a power of 2
multiplied by a mantissa that is between 1 and 2,
like this:
\begin{verbatim}
2^(-2) * 1.0101010101010101010101[0,1].
\end{verbatim}
This shows the single-length format,
which has 23 bits after the binary point.

Long strings of bits
are usually shown in hexadecimal notation,
which uses the 16 characters
\begin{verbatim}
0, 1, 2, 3, 4, 5, 6, 7, 8, 9, a, b, c, d, e, f
\end{verbatim}
to represent four consecutive bits at a time.
As there are only 23 bits to be displayed, we represent the first
three bits in octal notation,
with the characters
\begin{verbatim}
0, 1, 2, 3, 4, 5, 6, 7.
\end{verbatim}
Thus the 23 bits after the binary point
are shown as an octal character
followed by five hexadecimal characters.
Thus, the narrowest interval containing $1/3$ is
\begin{verbatim}
2^(-2) * 1.2aaaa[a,b].
\end{verbatim}
In this way we get the table in Figure~\ref{fig:sioresults}.
\begin{figure}[htbp]
\begin{verbatim}
1/i     floating-point number     narrowest interval
        produced by standard      containing 1/i
        I/O library via compiler
--------------------------------------------------------
1/2     2^(-1) * 1.000000         2^(-1) * 1.000000[,]
1/3     2^(-2) * 1.2aaaab         2^(-2) * 1.2aaaa[a,b]
1/4     2^(-2) * 1.000000         2^(-2) * 1.000000[,]
1/5     2^(-3) * 1.4ccccd         2^(-3) * 1.4cccc[c,d]
1/6     2^(-3) * 1.2aaaab         2^(-3) * 1.2aaaa[a,b]
1/7     2^(-3) * 1.124925         2^(-3) * 1.12492[4,5]
1/8     2^(-3) * 1.000000         2^(-3) * 1.000000[,]
1/9     2^(-4) * 1.638e39         2^(-4) * 1.638e3[8,9]
1/10    2^(-4) * 1.4ccccd         2^(-4) * 1.4cccc[c,d]
1/11    2^(-4) * 1.3a2e8c         2^(-4) * 1.3a2e8[c,d]
\end{verbatim}
\caption{
Why standard input cannot be relied on to obtain an interval
for some of the common fractions.
}
\label{fig:sioresults}
\end{figure}

In this paper we develop algorithms to support software
that performs correct input and output of intervals.

\section{Previous work}

The problems described in the introduction have motivated
Rump \cite{rmp99} to develop a method for decimal input and output
for intervals. He computes two 2-dimensional arrays
{\tt low} and
{\tt upp}
of double-length floating-point numbers that satisfy
$$
\mbox{low}[d,e]
\leq
d \times 10^e
\leq
\mbox{upp}[d,e]
$$
for $d \in \{1,2,\dots,9\}$ and
for $e \in \{-340,-339,\dots,308\}$.
These two arrays contain a total of 11682 double-length floating-point numbers.

With these precomputed numbers available, a decimal numeral
$0.d_1 d_2 \dots d_k \times 10^e$
is converted on input to the interval
\begin{equation}
[ \sum_{i=1}^{k}\mbox{low}[d_i,e-i], \sum_{i=1}^{k}\mbox{upp}[d_i,e-i] ]
\end{equation}
of double-length floating-point numbers,
where the inner brackets enclose array subscripts.
The additions for the lower (upper) bound are performed in
downward (upward) rounding mode.
As a result of the possibly occurring roundings,
the interval obtained is not in general as narrow as possible.
To minimize the inevitable widening, Rump recommends performing
the  additions starting with the smallest array elements.

The work of Rump should be compared and contrasted with what we
will call here \emph{general-purpose} conversion algorithms.
Though these are not specifically intended for interval conversions,
they can perhaps be adapted to this purpose because of their guarantees on
accuracy.

Let us briefly review these algorithms.
Steele and White \cite{stlwht90}
formulated \emph{identity requirements} for conversions between
internal and external numerals.
The internal identity requirement is that conversion of an internal numeral
to external and back results in the same.
If the external numerals are unlimited in length,
this requirement can always be met.
Similarly,
the external identity requirement is that conversion of an external numeral
to internal and back results in the same.
Usually, the internal numerals have a fixed length,
so this requirement can only be met when the external numeral is not
longer than is warranted by this fixed length.

Every binary numeral is representable as a decimal one.
Hence,
if the internal numerals are binary and the external ones are decimal,
then it is trivial to satisfy the internal identity requirement by using
the decimal equivalent of the internal binary numeral.
The problem addressed by Steele and White is to do so with as few
decimal digits as possible.

To transfer binary floating-point numerals from one computer to
another by means of binary files, one needs to be assured that 
these files are written and read in a compatible way.
Moreover, the layout of the floating-point numerals of both machines
needs to be the same. Adherence to the IEEE floating-point format
is not enough: in addition, both machines have to be big-endian or
both need to be little-endian.
Because of these difficulties it is attractive
to convert internal binary numerals
to an external text file containing decimal numerals
and to have an algorithm
that guarantees faithful translation back to internal format.

Such use of decimal I/O requires that the conversion is efficient.
For this reason, Gay \cite{gay90correctly} and Burger and Dybvig \cite{burger96printing}
devised a faster decimal output.
The identity requirement assumes sufficiently accurate decimal input.
This was addressed by Clinger \cite{clinger90how} and by Gay \cite{gay90correctly}.

Compared to Rump's approach,
this has the advantage of not requiring a database
of precomputed numbers.

The identity requirements satisfied by
\cite{stlwht90,gay90correctly,burger96printing,clinger90how}
are convincing for the purposes envisaged by these authors.
However, the requirements for decimal interval I/O
are equally compelling and quite different:
\begin{enumerate}
\item
On input, compute the narrowest floating-point interval containing the decimal
numeral.
\item
On output, compute the narrowest interval containing the floating-point number
that is bounded by decimal numerals of specified length.
\end{enumerate}
To satisfy these requirements it seemed simplest to develop our algorithms
from first principles rather than to attempt to adapt the work referenced
above.

From \cite{hhkr93} it seems that XSC does the conversions
between binary and decimal numerals correctly, but the authors
do not say how it was done. Although the source code of XSC
is publicly available (GPL license), 
it contains more than $150$ files.
Moreover, the code is not documented
in such a way as to facilitate finding
the pieces of code that do the conversion, to gather them, 
and to understand how the conversion was done.

Our task, then, is to explain how to do the conversion,
and to show how to implement it.
As the utility of this kind of work requires executable programs,
we needed to pick an implementation language.
The primary language for numerical work is Fortran.
However, this is more in the direction of systems programming,
for which C/C++ is a reasonable choice of language.

\section{Decimal input}

The ``scientific notation'' for a number
consisting of a decimal fraction and an exponent
is almost universally used.
Thus it would seem that one only has to cater to this format
for I/O routines.
Yet for input there is a strong case to be made for
the pre-scientific notation of a fraction as a pair of integers.
Therefore we consider these two in turn.

\subsection{Rational fractions}
\label{sec:ratFrac}
In many situations the most convenient way
to input a number is as a rational fraction $p/q$,
where $p$ and $q$ are integers.
Not only is it convenient,
but there is also a convincing correctness criterion:
as there exists, in a given format,
a unique least floating-point interval
that contains $p/q$,
the input function should yield this interval.

An algorithm for this purpose is one that has been widely,
but not universally, taught to children for at least two centuries.
We will illustrate this algorithm with $p/q = 3/7$.
As this number is less than $1$, we know that the binary
fraction has the form \verb+0+$\cdot$\verb+...+.
How do we get the missing digits?
Multiply by two to get $6/7$.
As the result is less than one,
the result has the form
\verb+0.0...+
Multiply by two to get $12/7$.
As the result is not less than one,
the result has the form
\verb+0.01...+
and subtract 1 from $12/7$, so that we continue with $5/7$.
Double again, get $10/7$, so that we continue with $3/7$
after noting that
the result has the form
\verb+0.011...+
We already know what comes after that, so we have determined that the binary
equivalent of $3/7$ is \verb+0. 011 011 011 ...+

Let us check our computation. The first group of digits after the decimal
point is worth $3/8$.
Every next one is worth $1/8$ times the previous group of three.
So the value of the infinite string is
$(3/8)*(1 + (1/8) + (1/8)^2 + \ldots)$.
Let $S = 1 + (1/8) + (1/8)^2 + \ldots$.
Then
$S/8 = (1/8) + (1/8)^2 + \ldots = S-1$.
Hence $S = 8/7$ and we find that the infinite string is
$(3/8)*(8/7) = 3/7$.

Let us now translate this procedure to a machine-executable algorithm.
We can decide whether the next binary digit should be a 0 or a 1
by computing in a variable \verb+pwr+ (from ``power'')
the successive powers $2^{-k}$ for $k = 0, 1, 2, \ldots$
and adding some of these powers in a variable named \verb+frac+
(from ``fraction'').
We compute the next value of \verb+pwr+ at every step,
but only add it to \verb+frac+ when deciding to write a \verb+1+
in the algorithm described above.
This idea is embodied in the following code.

\begin{verbatim}
  float frac = 0.0;
  // fraction to be built up out of powers of two
  float pwr = 1.0; // power of 2 to add to fraction
  // let r = p/q
  while (frac + pwr > frac && p > 0) {
    // p > 0 and r - frac = (p/q)*pwr
    pwr = pwr/2.0; p = 2*p;
    if (p >= q) {
      p = p-q; frac = frac+pwr;
    }
  // if p = 0 then r = frac
  // if p > 0 then 0 < r - frac <= pwr
  }
\end{verbatim}
Typically, there are infinitely many binary digits.
Only the first of a finite segment can be accommodated in a floating-point
number.
At every iteration \verb+pwr+ is halved.
At some point this quantity becomes insignificant.
The criterion for this is whether adding \verb+pwr+ to \verb+frac+
makes any difference to \verb+frac+.
As soon as that is not the case,
the iteration terminates.

The last assertion states
that every time around the loop
we have that \verb+frac+ is a lower bound for $p/q$
and that the difference between the two
is at most $2^{-i}$,
where $i$ is the number of times around the loop.
This code builds in \verb+frac+ a lower bound to
$p/q$ that approaches $p/q$ as closely as the precision allows.
If $p/q$ has can be represented in the floating-point number format,
this shows by \verb+p+ becoming 0.

The above segment of code is the kernel of the function shown
in Figure~\ref{fig:pqCpp}
for the IEEE standard 754 single-length format.
It produces the least floating-point interval that contains
the fraction $p/q$.

\begin{figure}[htbp]
\begin{verbatim}
Interval* convert(int p, int q) {
// Assumes p and q are positive and less than 2^{30}.
// Returns the greatest float that is not greater than
// the rational r = p/q.
  float sf = scaleFactor(p,q);
  float frac = 0.0;
  // fraction to be built up out of powers of two
  float pwr = 1.0; // power of 2 to add to fraction
  while (frac + pwr > frac && p > 0) {
    // p > 0 and r - frac = (p/q)*pwr
    pwr = pwr/2.0; p = 2*p;
    if (p >= q) {
      p = p-q; frac = frac+pwr;
    } // p > 0 and r - frac = (p/q)*pwr
  }
  frac = frac*sf;
  if (p==0) // r = frac
    return new Interval(frac, frac);
  else return new Interval(frac, next(frac));
}
\end{verbatim}
\caption{\label{fig:pqCpp} 
A C++ function to convert a rational $p/q$ to a floating-point number.
The function
{\tt next()}
produces the least floating-point number
greater than its argument.
}
\end{figure}

The function assumes that $0.5 \leq p/q < 1$,
which is of course not in general the case.
It therefore needs the function in Figure~\ref{fig:scaleFactor}.

\begin{figure}[htbp]
\begin{verbatim}
float scaleFactor(int& p, int& q) {
// Let r = p/q.
// Returns a scale factor sf such that p/q is in [0.5,1)
// and r = (p/q)*sf
  float sf = 1.0; // the scale factor
  // r = (p/q)*sf; assume q < 2^{30} to avoid overflow
  while (q > p) { p = 2*p; sf /= 2.0; }
  // r = (p/q)*sf and q <= 2*p;
  // assume p < 2^{30} to avoid overflow
  while (p >= q) { q = 2*q; sf *= 2.0; }
  // r = (p/q)*sf and p < q <= 2*p;
  // therefore 0.5 <= (p/q) < 1
  return sf;
}
\end{verbatim}
\caption{\label{fig:scaleFactor} 
A C++ function to scale the fraction $p/q$.
}
\end{figure}

\subsection{Input of decimal floating-point numerals}
Many programming languages and data files
use a similar format for fractional numbers.
Though details may vary,
it is easy to extract from such files
an integer $e$ containing the exponent
and a string containing
the decimal digits $d_1, \ldots, d_n$ of the fraction,
where $d_1 \neq 0$.
We assume that these source code or data file elements are
intended to denote the rational number
$r = 10^e\sum_{i=1}^n d_i 10^{-i}$.
The function we aim at here takes as input $e$ and
$d_1, \ldots d_n$ and returns the narrowest
floating-point interval containing $r$.

Such numerals can be treated analogously to the rational fractions
discussed in Section~\ref{sec:ratFrac}.
A difference is that instead of
subjecting a rational fraction $p/q$ to repeated
doubling possibly combined with subtracting 1,
we do this with a string of decimal digits representing the mantissa of the
number to be input. 

A more significant difference
compared to the rational-fraction case
is the presence of a power of 10.
We need to convert to binary not only the mantissa, 
but also the input power of 10.
This happens in a preliminary stage
we call \emph{binarization}.

For example, suppose we desire to convert
$0.0123$ to binary.
According to our assumed convention,
we have $e = -1$
and have $d_1, \ldots d_n$ in the form of the string \verb+"123"+.

In this example the power of 10 is exchanged with power of 2 and
a corresponding change in mantissa as follows:
$
0.123 * 10^{-1} =
0.246 * 10^{-1} * 2^{-1} =
0.492 * 10^{-1} * 2^{-2} =
0.984 * 10^{-1} * 2^{-3} =
0.1968 * 2^{-4}.
$

To implement binarization,
we need a function to double a decimal mantissa given
as a string of decimal digits.
Each digit is doubled by integer arithmetic operating on
the numerical equivalent of the digit.
In this operation neither rounding nor overflow can occur.
See the function \verb+mul2+ in Figure~\ref{fig:mul2}.
This function is not a general-purpose doubling routine:
it is specific to its argument being the decimal digits
$d_1, \ldots, d_n$
of a mantissa of the form
$\sum_{i=1}^n d_i 10^{-i}$.
Accordingly, when the result is $1$ or greater,
this additional digit is not inserted into the resulting string.
Instead the last carry is returned.
By inspecting it, the calling code can determine whether the mantissa
has overflowed.

\begin{figure}[htbp]
\begin{verbatim}
int mul2(string& mnts) {
// Multiplies decimal mantissa in argument by 2.
// Returns last carry.
  int carry = 0;
  int dd; // Result of doubling a decimal digit.
  for(int i = mnts.length()-1; i >= 0; i--){
    dd = 2*(mnts[i]-'0') + carry;
    mnts[i] = dd%10 + '0'; carry = dd/10;
  }
  if (mnts[mnts.length()-1] == '0')
    mnts.erase(mnts.length()-1,1);
  // mnts is a mantissa, hence no trailing zeros
  return carry;
}
\end{verbatim}
\caption{\label{fig:mul2}
A C++ function to double a number of the form
$\sum_{i=1}^n d_i 10^{-i}$ where the $d_1,\ldots,d_n$
are given in the string argument {\tt mnts}.
}
\end{figure}

With the doubling function available,
the binarization function is straightforward.
See Figure~\ref{fig:binarize}.
\begin{figure}[htbp]
\begin{verbatim}
void binarizeExp(string& mnts, int& exp) {
// When called, exp is a decimal exponent.
// On exit, exp is a binary exponent and mnts is adjusted,
// so that the same number is denoted.
  int binExp = 0; int carry = 0;
  while (exp < 0 ) {
    carry = mul2(mnts); binExp--;
    if (carry != 0) { exp++; mnts.insert(0,"1"); }
  }
  while (exp > 0 ) {
    div2(mnts); binExp++;
    if (mnts[0] == '0') { exp--; mnts.erase(0,1); }
  }
  // exp = 0
  exp = binExp;
}
\end{verbatim}
\caption{\label{fig:binarize}
A C++ function to binarize a number of the form
$10^e\sum_{i=1}^n d_i 10^{-i}$ where the digits $d_1,\ldots,d_n$
are given in the string argument {\tt mnts}.
If we denote the values on exit by adding primes,
we have
$
10^e\sum_{i=1}^n d_i 10^{-i}=
2^{e'}\sum_{i=1}^{n'} {d'}_i 10^{-i}
$
}
\end{figure}

In the next stage, \emph{normalization}, we
ensure that the mantissa is in the interval $[0.5,1)$.
In our current example this happens by means of the steps
$
0.1968 * 2^{-4} =
0.3936 * 2^{-5} =
0.7872 * 2^{-6}
$
The function for normalizing is in Figure~\ref{fig:normalizeEM}.
\begin{figure}[htbp]
\begin{verbatim}
void normalize(string& mnts, int& exp) {
// Maintaining the value of the denoted number,
// adjusts exp so that 0.5 <= 0.mnts < 1.
  while ((mnts[0]-'0') < 5) {
    exp--; mul2(mnts); // discard zero carry
  }
}
\end{verbatim}
\caption{\label{fig:normalizeEM}
A C++ function to normalize a number of the form
$\sum_{i=1}^n d_i 10^{-i}$ where the mantissa $d_1,\ldots,d_n$
is given in the string argument {\tt mnts}.
}
\end{figure}

After binarization and normalization we are ready
to start the conversion of the decimal mantissa to binary.
A good starting point is one of the radix conversion methods given by
Knuth \cite{knth69}, section 4.4, where he converts
from radix $b$ to radix $B$.
In our case we have $b = 10$ and $B = 2$.
The options are
to divide by 10 using radix-2 arithmetic
and
to multiply by $B = 2$ using radix-$10$ arithmetic.
We select the latter.

Thus we have a fractional number $u$ given as a string of decimal digits.
Knuth obtains the digits $U_1, U_2, \ldots$ of the binary representation as follows:
\begin{eqnarray*}
U_1   &=& \lfloor uB \rfloor             \\
U_2   &=& \lfloor \{uB\}B \rfloor        \\
U_3   &=& \lfloor \{\{u\}B\}B \rfloor    \\
\ldots,                                   \\
\end{eqnarray*}
where $\{x\}$ denotes $x \bmod 1$,
which is  $x - \lfloor x \rfloor$.
 
We must not only obtain the successive binary digits
$U_1$,
$U_2$, and
$U_3$, but we must also pack them as a floating-point number.
Hence we modify Knuth's iteration to the following,
of which we show the first few steps,
continuing the previous example.
\begin{eqnarray*}
0.7872 & = & 0.7872 * 1           \\
       & = & 1.5744 * (1/2)           \\
       & = & 1/2 + 0.5744 * (1/2)      \\
       & = & 1/2 + 1.1488 * (1/4)      \\
       & = & 1/2 + 1/4 + 0.1488 * (1/4)      \\
       & = & 1/2 + 1/4 + 0.2976 * (1/8)      \\
       & = & 1/2 + 1/4 + 0.5952 * (1/16)      \\
       & = & 1/2 + 1/4 + 1.1904 * (1/32)      \\
       & = & 1/2 + 1/4 + 1/32 + 0.1904 * (1/32),      \\
\ldots \\
\end{eqnarray*}
Giving \verb+0.11001...+ as
the binary representation of $0.7872$.
We see that in this way a sum of powers of 2 is built up while the remaining
decimal mantissa is multiplied by an ever smaller factor.
The iteration is terminated when this product is less than the machine
precision.
At that point the sum of the powers of 2 is the left bound of the narrowest
interval containing 0.7872.
The multiplications and additions, though floating-point operations,
are, by their special nature,
performed without rounding error.

In the final stage,
we scale the lower bound by the factor $2^{-6}$
resulting from binarization and normalization.
We perform this scaling by successive divisions or multiplications by 2,
again assuring the absence of rounding errors.
Figure~\ref{fig:decInp} displays a function along these lines that finds the
narrowest floating-point interval.
\begin{figure}[htbp]
\begin{verbatim}
interval convertFrac(string& mnts, int& exp) {
  interval result; // interval to be returned
// r == 10^{exp} * 0.mnts
  binarizeExp(mnts, exp);
// r == 2^{exp} * 0.mnts
  normalize(mnts, exp);
// r == 2^{exp} * 0.mnts and 0.5 <= 0.mnts < 1
// r' == 0.mnts and 0.5 <= 0.mnts < 1
  float frac = 0.0; float pwr = 1.0;
  int carry = 0;
  while (frac+pwr > frac && mnts.length() > 0) {
  // 0 <= r' - frac = 0.mnts * pwr
    pwr /= 2.0; // pwr == 2^{-i}
    carry = mul2(mnts);
    if (carry != 0) { frac += pwr;
    // subtract 1 by discarding nonzero carry
    }
  } // 0 <= r' - frac = 0.mnts * pwr
    // and pwr is negligible compared to frac
    // hence frac is the greatest flpt less than r'
  // scale frac with exp
  while (exp > 0) {frac *= 2.0; exp--;}
  while (exp < 0) {frac /= 2.0; exp++;}
  result.lb = frac; 
  result.ub = (mnts.length() == 0) ? frac : next(frac);
  return result;
}
\end{verbatim}
\caption{\label{fig:decInp}
A C++ function to find the narrowest floating-point interval that contains
a given number of the form
$10^e\sum_{i=1}^n d_i 10^{-i}$ where the $d_1,\ldots,d_n$
are given in the string argument {\tt mnts} and $e$ is given in the argument
{\tt exp}.
}
\end{figure}
This function that the fraction $x$ to be converted
is non-negative.
In that case the result is $[a,a']$, where $a'$ is the least
floating-point number greater than floating-point number $a$.
This function can also be used for a negative fraction $y$.
If our function gives $[a,a']$ with input $-y$,
then we change the output to $[-a',-a]$.

This function uses some auxiliary declarations:
one that defines \verb+interval+
and one that defines a function to determine
the next floating-point number after a given one.
The auxiliary definitions are displayed in Figure~\ref{fig:aux}.
\begin{figure}[htbp]
\begin{verbatim}
typedef struct interval { float lb, ub; };

float next(float x) {
// returns the least flpt number greater than nonnegative x
// if x is normalized and if it is less than the greatest
// flpt number
  // x0 = x
  float sf = 1.0; // becomes the scale factor
  // x0 = x*sf
  while (x < 1.0) { x *= 2.0; sf /= 2.0; }
  // x0 = x*sf and 1 <= x
  while (x >= 2.0) { x /= 2.0; sf *= 2.0; }
  // x0 = x*sf and 1 <= x < 2
  float eps = 1.0; // becomes the machine epsilon
  while (((float)1.0 + eps) > (float)1.0) eps *= 0.5;
  // eps is the first one that did not make a difference
  eps*=2.0;
  // eps is the last one that did make a difference
  // by definition the machine epsilon
  return (x+eps)*sf;
}
\end{verbatim}
\caption{\label{fig:aux}
Some definitions auxiliary to Program~\ref{fig:decInp}.
}
\end{figure}

\section{Decimal output}

As in decimal input,
the main concern is to avoid rounding errors.
This is of course taken care of by representing all digits,
binary or decimal,
separately as small integers.
Operations on these are free from rounding errors
because the operands are integer;
they are immune to overflow because of their smallness.
However, we would also like to speed up conversion
as much as possible by using operations on floating-point numbers
when we can be sure that no rounding errors occur
or by using operations of full-size integers
when we can be sure that no overflow can occur.

Our starting point is the assumption
that a floating-point number can be represented
by an integer for the exponent part
and by an integer that is represented
by the same sequence of bits as there are
in the mantissa of the floating-point number.
In the case of the IEEE standard single-length format,
this integer is 24 bits with the most significant bit equal to 1.
That is, for given floating-point $f$,
we find integers $e$ and $m$ such that $f = m*2^e$.
As before, we use the case of the IEEE standard 754
single-length floating-point format as example.

Consider the assignment \verb+m = f+, which is legal in C/C++.
It implicitly converts the floating-point number \verb+f+ to an
integer if \verb+m+ is an integer.
If \verb+f+ is not an integer, then this assignment does not result in
\verb+m+ containing the bits of the mantissa of \verb+f+.
If \verb+f+ is an integer, then this assignment may also not result in
\verb+m+ containing the bits of the mantissa of \verb+f+:
\verb+f+ may be larger than the largest integer.
However, if we ensure that \verb+f+ is in the interval
$[2^{23}, 2^{24})$, it is both assured to be an integer,
while at the same time not causing overflow when assigned to a 32-bit
integer.
Thus we may extract the bits of \verb+f+ and place them in \verb+m+
by simply performing the assignment \verb+m = f+
provided that we first scale \verb+f+ to be in this range.
The scaling is performed by doubling or halving of \verb+f+
and therefore does not cause rounding errors.
The required scaling operations are accounted for in the binary
exponent \verb+e+.
See Figure~\ref{fig:fltoem} for a function to perform the conversion
from $f$ to the corresponding $e$ and $m$.

\begin{figure}[htbp]
\begin{verbatim}
void fltoem(const float& f0, int& e, int& m) {
// returns m and e such that f0 = m * 2^e
// with 2^{23} <= m < 2^{24} that is
// m has the same bits as 1.significand of f0
  float f = f0;
  float e23 = 1.0;
  for (int i = 0; i < 23; i++) e23 *= 2.0; // e23 = 2^{23}
  float e24 = 2.0 * e23;
  e = 0; // f0 == f * 2^e
  while (f <  e23) {f *= 2.0; e--;} // f0 == f * 2^e
  while (f >= e24) {f /= 2.0; e++;} // f0 == f * 2^e
  // f is an integer and 2^{23} <= f < 2^{24}
  m = (int)f; // f0 == m*2^e and 2^{23} <= m < 2^{24}
}
\end{verbatim}
\caption{\label{fig:fltoem}
A function to realize the relation $f_0 = m*2^e$,
with $f_0$ as input and $m$ and $e$ as output.
Example for single-length floating-point format.
}
\end{figure}

It remains to determine the string of decimal digits that is equal to
$m*2^e$ with $m$ and $e$ given as variables of type integer.
If we place no limit on the length of the string,
this is always possible: decimal strings are closed under doubling and halving.
We first convert the integer variable $m$ to a decimal string.

This is done by one of the best known algorithms in programming.
One of the forms in which it appears is the function \verb+itoa+
in \cite{krnrtch88}. A streamlined version is the following,
though it prints the digits in reverse order:
\begin{verbatim}
void print(int n) {
  while (n>0) { cout << n%10; n = n/10; }
  cout << endl;
}
\end{verbatim}

What we need is basically the same algorithm,
but elaborated by the need
to output to a string rather than to print.
We also need to ensure
that the digits are not reversed.
Thus we arrive at the function in Figure~\ref{fig:decimalizeMnts}.
\begin{figure}[htbp]
\begin{verbatim}
void decimalizeMnts(int& n, string& s) {
  s.erase(0, s.length()); // ensures s is empty
  while (n > 0) {
    s.insert(s.begin(),1,(char)('0' + (n%10)));
    // inserts digit n%10 at the beginning of s
    n = n/10;
  }
}
\end{verbatim}
\caption{\label{fig:decimalizeMnts}
Decimalizes the integer mantissa.
}
\end{figure}

The result of \verb+decimalizeMnts+ is the representation of the
\emph{integer} \verb+m+ such that $ f = m*2^e$.
What we want is a \emph{mantissa}.
That is,
we need to move the decimal point from the right to the left.
This is done by multiplying by a suitable power of 10,
which is the length of the decimal numeral.
Hence the function in Figure~\ref{fig:normalize}.

\begin{figure}[htbp]
\begin{verbatim}
int normalize(const string& mnts) {
// Returns the decimal exponent necessary to turn mnts
// from a whole number into a normalized mantissa.
// Assumes mnts has no leading zeros.
  return mnts.length();
}
\end{verbatim}
\caption{\label{fig:normalize}
Transforms the integer mantissa to a fractional one.
}
\end{figure}

In the final stage we convert the binary exponent \verb+e+
to increments or decrements of the decimal exponent,
with a concomitant change in the mantissa.
This happens in the function shown in Figure~\ref{fig:decimalizeExp}.
\begin{figure}[htbp]
\begin{verbatim}
void decimalizeExp(string& mnts, int& binExp, int& decExp) {
// Reduces the binary exponent binExp to zero.
// Assumes that mnts has no leading zeros.
  while (binExp > 0) {
    binExp--;
    if (mul2Mnts(mnts) == 1) {
      decExp++; mnts.insert(0,"1");
    }
    if (mnts[mnts.length()-1] == '0')
      mnts.erase(mnts.length()-1,1);
  }
  while (binExp < 0) {
    div2(mnts); binExp++;
    if (mnts[0] == '0') {
      mnts.erase(0,1); decExp--;
    }
  }
}
\end{verbatim}
\caption{\label{fig:decimalizeExp}
Decimalizes the binary exponent.
}
\end{figure}

To summarize, we transform a floating-point number $f$ to an integer
$exp$ and a mantissa $mnts$ such that
$f = 0.mnts * 10^{exp}$
by first transforming it to an integer mantissa with a binary exponent
(function \verb+fltoem+),
then decimalizing the mantissa
(function \verb+decimalizeMnts+),
then transforming the integer mantissa to a fractional one
(function \verb+normalize+),
and finally folding the binary exponent into the decimal
exponent produced by the normalization stage.
This completes the decimal output of a floating-point number.
See the function \verb+decOut+ in Figure~\ref{fig:decOut}.

\begin{figure}[htbp]
\begin{verbatim}
void decOut(float f, int& exp, string& mnts) {
// returns the decimal string equal to f
  int e, m;
  fltoem(f, e, m); // f = m * 2^e
  decimalizeMnts(m, mnts);
  // mnts is decimal equivalent of m
  exp = normalize(mnts);
  decimalizeExp(mnts, e, exp);
}
\end{verbatim}
\caption{\label{fig:decOut}
Finds the decimal equivalent of floating-point number.
}
\end{figure}

The function \verb+decOut+ in Figure~\ref{fig:decOut}
gives all of the decimals of a numeral
that is equal to the floating-point number in question.
To obtain the lower bound of the required narrowest interval
bounded by decimal numerals of specified length,
the output of \verb+decOut+ needs to be truncated to this length.
The upper bound is then obtained
by adding to it one unit of the last decimal place.

\section{Conclusions}

To input an interval
given as a pair of decimal numerals,
one needs to assure
that the resulting pair of floating-point numbers
contains the input pair.
The fact that a decimal numeral
is typically not representable
as a binary floating-point number makes it necessary
that the internal result is wider than the input.
We want it to be wider by as little as possible.
The standard library functions
are not adequate for this purpose.
In this paper we have developed functions
that support such interval input.

For output, there is in principle no problem:
for every finite binary floating-point numeral
there is a decimal numeral
that is equal to it.
Therefore, if one is willing to accept lengthy numerals,
it is possible to obtain as output an interval
that is equal to the internally stored floating-point interval.
But to obtain an output interval
that is both correct and not too long,
it is in general necessary
to round the decimal bounds in the correct directions.

In this paper we have developed functions
that can be used as basic building blocks
for the required input and output for intervals.
These include
\begin{enumerate}
\item
For every decimal numeral,
computing the unique narrowest floating-point interval containing it.
\item
For every finite floating-point number,
the decimal numeral that is equal to it.
\end{enumerate}

In designing these building blocks,
we have been guided by the criteria of efficiency,
language-independence, and machine-independence.

\paragraph{Efficiency}
To avoid rounding errors and overflow
one can perform all arithmetic
on numbers represented by strings of digits.
These digits are small integers
and are therefore operated upon
without rounding error and without overflow.
The disadvantage is the resulting slowness of the algorithms.
We have therefore taken advantage of several situations
in which floating-point operations do not give rounding errors
and are much faster.

\paragraph{Language-independence}
The most extreme form of language-independence
is to write the algorithms in pseudo-code.
To ensure that our algorithms can be verified,
we have chosen a simple subset of C++
that almost fits inside the C language.
The few excursions beyond the bounds of C
are amply rewarded, for example,
by the availability of the \verb+string+ class.
Of course, the most important part of verification
is the understanding of the code.
In this respect C is not better than pseudo-code,
but not worse either.
The reason for presenting our algorithms in executable form
is that execution is a useful check
on verification by reading the code.

We have facilitated reading the code
by including as
assertions comments.
Of each of these we claim
that it holds whenever execution passes the assertion.

But the fact that we want executable algorithms,
and hence rely on a specific programming language,
does not imply that we feel free
to use all the facilities of that language.
For example, to determine the next floating-point number
after a given one,
we can use the facilities of C
to access the bits of a floating-point number.
But such facilities vary widely among programming languages.
We therefore prefer to determine
the next floating point number
by means of arithmetic operations,
so that the resulting algorithm
is more readily transferable to another language.

\paragraph{Machine-independence}
For input it is desired to find
the smallest floating-point interval
containing the number denoted by the input numeral.
Usually such a requirement
is met by modifying the default rounding mode of the processor.
We have refrained from using this option,
as the specification of the programming language
does not include these operations.
As a result, the different implementations of C
come with implementation-specific functions
to perform these operations.

The motivation of this paper was the usual situation
where the semantics of an interval is
that it contains one or a finite number of solutions.
In this situation correctness of the conversion process
means obtaining an interval
that \emph{contains} the original one.
It also happens that the semantics of the original interval
is that \emph{all} points in it are solutions,
as can happen with inequalities.
In such a situation the correctness criterion
is to transform the original interval
to one that is contained in it,
again differing by as little as possible.
The functions presented here
can serve as building blocks in this situation also.

\section*{Acknowledgements}
This research was supported by the University of Victoria and by
the Natural Science and Engineering Research Council of Canada.


\end{document}